\documentclass[12pt,reqno]{article}
\def\hybrid{\topmargin 0pt      \oddsidemargin 0pt
        \headheight 0pt \headsep 0pt
        \textwidth 135true mm       
        \textheight 224true mm         
        \marginparwidth 0.0in
        \parskip 0pt plus 1pt   \jot = 1.5ex}

\usepackage{amssymb}
\usepackage{amsmath}
\usepackage{amsthm}
\hybrid
\usepackage{amssymb}
\usepackage{amsmath}
\usepackage{amsthm}

\newcommand{\dst}{\displaystyle}

\newcommand{\lr}{L_{\overline{r}}}
\newcommand{\lp}{L_{\overline{p}}}
\newcommand{\lpr}{L^{\overline{\alpha}}_{\overline{p},\overline{r}}}

\newcommand{\ta}{T^{\overline{\alpha }}}
\newcommand{\sa}{S_{\overline{\alpha }}}

\newtheorem{remark}{Remark}[section]
\newtheorem{lemma}{Lemma}[section]
\newtheorem{theorem}{Theorem}[section]
\newtheorem{corol}{Corollary}[section]

\title{Characterization of anysotropic Liouville-type mixed spaces and its application}
\author{V. A. Nogin and E. E. Ournycheva}
\date{}

\begin{document}
\thispagestyle{empty}
\maketitle

\begin{abstract}
Within the framework of the method of hypersingular integrals we
obtain a characterization of the anysotropic spaces $\lpr$.
These spaces are defined to consist of functions $f(x)\in
L_{\overline{r}}$ for which
\[
 F^{-1}\left(\sum_{j=1}^n |\xi_j|^{\alpha_j}\right)Ff\in
 L_{\overline{p}},\quad \alpha_j>0.
\]
The mentioned characterization is applied to prove the
denseness of the class $C_0^\infty$ in $\lpr$.
\end{abstract}


\setcounter{equation}{0}
\setcounter{theorem}{0}
\setcounter{lemma}{0}
\section{Introduction}

We consider the anysotropic Liouville-type mixed spaces,
which are defined via Fourier transform as follows:
\begin{equation}
  \lpr\equiv L^{\overline{\alpha}}_{\overline{p},\overline{r}}({\Bbb R}^n)=
\left\{f: f\in L_{\overline{r}}, F^{-1}
\left(\sum_{j=1}^{n}|\xi_j|^{\alpha_j}\right)Ff\in
L_{\overline{p}} \right\},
\label{1}
\end{equation}
where $\alpha _i>0$, $\overline{p}=(p_1,\dots,p_n)$,
$\overline{r}=(r_1,\dots,r_n)$, $1\leq p_i,r_i<\infty $
($i=1,\dots,n$), $\lp\equiv L_{\overline{p}}({\Bbb R}^n)$ being
the well-known space equipped with the mixed norm
\begin{equation}
 \|f\|_{\overline{p}} =
 \left(\,\int\limits_{{\Bbb R}^1}\left(\,\int\limits_{{\Bbb
 R}^1}\dots \left(\,\int\limits_{{\Bbb
 R}^1}|f(x)|^{p_1}dx_1\right)^ {\frac{p_2}{p_1}} \dots
 dx_{n-1}\right)^{\frac{p_{n}}{p_{n-1}}}dx_n\right)^{\frac{1}{p_n}}.
\label{2}
\end{equation}
(We note that some basic properties of the spaces $\lp$ were
established in \cite{BenedekPanzone,BesovIlinNIk}.) The Fourier transform
$Ff$ in (\ref{1}) is treated in the sense of distributions on a suitable space of
test functions (see Section~\ref{new:sect2}).

In the isotropic case when $\alpha_1=\dots=\alpha_n=\alpha$,
$p_1=\dots=p_n=p$, and $r_1=\dots=r_n=r$ the spaces (\ref{1})
coincide with the spaces $L^\alpha_{p,r}$ of functions $f\in L_r$
having their Riesz derivatives ${\Bbb D}^\alpha f\in L_p$. Such spaces
were first introduced and investigated by S.\,G.~Samko in
\cite{b17,StudMath} (see also the books
\cite{Samko16} and \cite{Samko17}).
We note that
the spaces (\ref{1}) were treated by P.\,I.~Lizorkin
in the case $\overline{p}=\overline{r}$ (see \cite{MIAN}--\cite{b12}).
The spaces $L^{\overline{\alpha }}_{p,r}$ were studied by
A.\,A.~Davtyan in \cite{b3,b4} in the case $1<p<n/\alpha^*$,
$1/\alpha^*=\frac{1}{n}\sum\limits_{j=1}^n1/\alpha_j$, $1<r<\infty$, and by V.\,A.~Nogin
and G.\,P.~Emgusheva (see \cite{b6,b7}) for $1\leq p,r<\infty$.

Here we deal with the most general anysotropic case of vector-valued
$\overline{\alpha}$, $\overline{p}$, and $\overline{r}$.
We obtain a characterization of the spaces $\lpr$
via anysotropic hypersingular integrals (HSI)
\begin{equation}\label{new:1.3}
 (\ta f)(x)=
 \lim_{\stackrel{\varepsilon\to0}{(L_{\overline{p}})}}
 (T_\varepsilon^{\overline{\alpha}}f)(x),
\end{equation}
introduced by P.\,I.~Lizorkin (see \cite{b11}),
where $T_\varepsilon^{\overline{\alpha}}f$ is the "truncated"
integral
\begin{equation}\label{new:1.4}
  (\ta_\varepsilon f)(x)=
 \int\limits_{\rho(t)>\varepsilon}\frac{(\Delta_t^{2\ell}f)(x)}
{\rho ^{n+\alpha ^*}(t)}\,dt,\quad 2\ell>\max_j \alpha_j.
\end{equation}
Here $(\Delta_t^\ell f)(x)$ is the centered  finite difference of the
function $f(x)$.
The function $\rho(t)$, which is referred to as "anysotropic distance",
is a positive solution of the equation
$\sum\limits_{i=1}^n x_i^2\rho^{-2\lambda_i}=1$,
$\lambda_i=\alpha^*/\alpha_i$.
Namely, we prove that
\begin{equation}\label{new:1.5}
  \lpr=\{f: f\in \lr,\ \ T^{\overline{\alpha}}f\in \lp\},
\end{equation}
$\alpha_i>0$, $1<p_i<\infty$, $1\leq r_i<\infty$, under some additional
 restrictions on $\overline{\alpha}$ in the case $\alpha^*\geq n$ (see (\ref{2.6})).

We also give an application of the mentioned characterization.
Basing ourselves on equality
%With the aid of
(\ref{new:1.5}), we prove that the class
$C_0^\infty$ is dense in $\lpr$, $\alpha_i>0$, $1<p_i<\infty$,
$1\leq r_i<\infty$ (under the mentioned restrictions on $\overline{\alpha}$ in the case
$\alpha^*\geq n$). It should be noted that the denseness of the
class $C_0^\infty$ in $L_{p,r}^\alpha$ was proved in \cite{b17}
for $1<p<n/\alpha$ and in \cite{b13} for $p\geq n/\alpha$. In the
anysotropic case the corresponding result was established
in \cite{DGTU98}, where the authors assumed that
$r_i>p_i$ or $r_i\leq p_i$ for every $i$, $1\leq i\leq n$, if
$\sum\limits_{i=1}^n\frac{1}{\alpha_ip_i}\leq 1$.
One of the goals of this paper is to get rid of these unnatural restrictions.

We note that the proof of denseness given in Theorem~\ref{new:t4.1} looks very
non-trivial. We first prove the statement of this theorem in the
case $0<\max\limits_j \alpha_j<1$. After that, we extend it to the
remaining values of parameters with the aid of induction.

We also point out the papers \cite{b23} and \cite{b24}, where the author dealt
with some generalization of the spaces of
$L_{p,r}^{\overline{\alpha}}$-type to the case of a more general
anysotropic distance (in comparison with that used in (\ref{new:1.4})).

The paper is organized as follows. In Section~\ref{new:sect2} we
formulate our main results, see
Theorems~\ref{new:t3.1}  and ~\ref{new:t4.1}.
Sections~\ref{sect1} and  \ref{sect2}
can be regarded as a background to the proof of the mentioned theorems.
They contain the necessary preliminaries and some auxiliary statements
respectively. In Section~\ref{new:sect5} we prove Theorems ~\ref{new:t3.1}  and ~\ref{new:t4.1}.
For the sake of convenience we gather some technical fragments of
the proofs in Appendix. Some results presented in this paper were announced in
\cite{DGTU98}.

\setcounter{equation}{0}
\setcounter{theorem}{0}
\setcounter{lemma}{0}
\section{The main results}\label{new:sect2}

We define the Liouville-type spaces  $\lpr$, $1\leq r_i<\infty $,
$1\leq p_i<\infty $, $\alpha _i>0$, $i=1,\dots ,n$, by
equality (\ref{1}), where the Fourier transform $Ff$ is interpreted in
the sense of  $\Psi'$-distributions (see Subsection $3.1$). We put
\[
 \| f\|_{\lpr}=\| f\|_{\overline{r}}+
 \left\|F^{-1}\left(\sum\limits_{j=1}^{n}
 |\xi _j|^{\alpha _j}\right)
 Ff\right\|_{\overline{p}}.
\]
The following theorem provides a characterization of the space $\lpr$ by means of
HSI (\ref{new:1.3}).

\begin{theorem}\label{new:t3.1}
Let $1<p_i<\infty$, $1\leq r_i<\infty$, $\alpha_i>0$, $i=1,\dots,n$, and
\begin{equation}\label{2.6}
\begin{array}{c}
 \dst \sum\limits_{j=1}^{n}\frac{1+k_j}{\alpha _j}\ne 1,\quad
  |k|=0,1,\dots ,m-1,\ \  m=[\gamma ], \\
 \dst \gamma=\max\limits_{j}\alpha_j
 \left(1-\sum\limits_{j=1}^{n}\frac{1}{\alpha_j}\right).
\end{array}
\end{equation}
Then  equality (\ref{new:1.5}) holds. Moreover,
the norms $\displaystyle \|f\|_{\lpr}$ and
$\|f\|_{\overline{r}}+\|T^{\overline{\alpha}}f\|_{\overline{p}}$
are equivalent.
\end{theorem}
\begin{remark}\label{r2.1}
For the rest of the paper we assume the condition (\ref{2.6}) to be fulfilled.
\end{remark}
We apply the characterization (\ref{new:1.5}) to  prove the next theorem.

\begin{theorem}\label{new:t4.1}
Let $1<p_i<\infty $, $1\leq r_i<\infty $, $\alpha _i>0$,
$i=1,\dots ,n$.
Then the space $C_0^\infty $ is dense in $\lpr$.
\end{theorem}

\setcounter{equation}{0}
\setcounter{theorem}{0}
\setcounter{lemma}{0}

\section{Preliminaries}\label{sect1}
\subsection{Notation:} $S$ is the Schwartz space of
rapidly decreasing smooth functions; $C_0^\infty $ is its subspace
of finite functions; $\Psi_0$, $\Psi $ are the Lizorkin spaces of
functions in $S$, vanishing together with all their derivatives at the
origin and on the coordinate hyperplanes respectively; $\Phi _0$,
$\Phi$ are their Fourier duals;
$\Psi'$, $\Phi _0'$ are the corresponding spaces of distributions;
$C_0^\infty(0,\infty)$ is
the class of finite $C^\infty$-functions on the semi-axis
$(0,\infty)$  supported outside zero; $M_{\overline{p}}$ is the class of $\overline{p}$-
multipliers; $X \rightarrow Y$ is the continuous embedding of the normed space $X$
into the normed space $Y$;
$(\Delta^\ell_yf)(x)=\sum\limits_{k=0}^
{\ell}(-1)^kC_\ell^kf\bigl( x+(\ell/2-k) y\bigr)$ is the centered
finite difference of the function $f(x)$ of order $\ell$ with the step $y$;
 $(\widetilde{\Delta}^\ell_yf)(x)=
\sum\limits_{k=0}^{\ell}(-1)^kC_\ell^kf(x-ky)$
is the non-centered finite difference; $\displaystyle (
Ff)(y)=\widehat{f}(y)=\int_{{\Bbb R}^n}f(x)e^{ix\cdot y}dx$ is
the Fourier transform of the  function $f(x)$; $\displaystyle
( F^{-1}f)(x)=\widetilde{f}(x)=\frac{1}{(2\pi)^n}\int_{{\Bbb
R}^n} f(y)e^{-ix\cdot y}dy$ is the inverse Fourier transform;
$\displaystyle \langle f,\omega\rangle=\int_{{\Bbb R}^n}
\overline{f(x)}\omega(x)dx$. All constants in various estimates
not necessarily the same at each occurrence are denoted by the
same letter $C$. The end of proof is denoted by $\blacksquare$.

\subsection{Estimates for some integrals and finite differences}
The following lemmas were proved in \cite{b17} in
the isotropic situation (when $\rho(x)=|x|$, $p_1=\dots=p_n$,
and $s_1=\dots=s_n$).
The proofs of their analogues
in the anysotropic case are much in lines with those given in \cite{b17}
and we omit them.

\begin{lemma}\label{l1.0}
Let $\gamma>0$, $1\leq s_i<\infty$, $i=1,\dots,n$,
\begin{displaymath}
I(t)=\left(\,
\int\limits_{{\Bbb R}^1}dx_n\dots
 \left(\,\int\limits_{{\Bbb
 R}^1}
\left(\prod\limits_{k=0}^\ell[1+\rho(x-kt)]^{-\gamma s_1}\right)
dx_1
 \right)^{s_2/s_1}\dots
\right)^{1/s_n},
\end{displaymath}
$\ell=1,2,\dots$. Then
\begin{equation}
I(t)\leq
C(1+\rho(t))^{\sum\limits_{i=1}^n\lambda_i/s_i-(\ell+1)\gamma},\quad if\
если\
\frac{1}{\ell+1}\sum_{i=1}^n\frac{\lambda_i}{s_i}<\gamma<
\sum_{i=1}^n\frac{\lambda_i}{s_i},
\label{1.1}
\end{equation}
\begin{equation}
I(t)\leq C(1+\rho(t))^{-\gamma\ell},\quad if \ если\
\gamma>\sum_{i=1}^n\frac{\lambda_i}{s_i}.
\label{1.2}
\end{equation}
\end{lemma}

\begin{lemma}\label{l1.3}
If $a(x)\in C_0^\infty$, then the following inequality holds:

\begin{equation}
\left|(\tilde\Delta_t^l a)(x)\right|\leqslant
C\Biggl[\frac{\rho^l(t)}{\prod\limits_{k=0}^l(1+\rho(x-kt))}\Biggr]^\theta,
\quad \ell=1,2,\dots,\  \theta=\min_i\left(\frac{\alpha^*}{\alpha_i}\right).
\label{1.5}
\end{equation}
$ x,t\in{\Bbb R}^n$.
\end{lemma}

\begin{lemma}\label{l1.4}
If $a(x)\in C_0^\infty$, then
\begin{equation}
\left\|\tilde\Delta_t^\ell a\right\|_{\overline{p}}\leq
C\left[\frac{\rho(t)}{1+\rho(t)}\right]^{\ell\theta},\quad
\ell=1,2,\dots,\quad \displaystyle
\theta=\min_i\left(\frac{\alpha^*}{\alpha_i}\right).
\label{1.6}
\end{equation}
\end{lemma}

\subsection{Some properties of the mixed spaces $\lp$}

We need the following properties of functions in $\lp$.

\begin{theorem}[\cite{BesovIlinNIk}]\label{t1.1}
Let $f(x)\in L_{\overline{p}}({\Bbb R}^n)$,
$1\leq p_i< \infty$, $i=1,\dots,n$,
%$f(x)\in L_{\overline{p}}({\Bbb R}^n)$
and let $a(t)$ be an averaging kernel, that is, $a(t)\in{L_{1}}$ and
$\displaystyle\int\limits_{{\Bbb R}^n}a(t)\,dt=1$. Then
\[
 \lim_{\delta\to0}\|f-f_\delta\|_{\overline{p}}=0,
\]
where
\begin{equation}
f_\delta(x)=\int\limits_{{\Bbb R}^n}a(t)f(x-\delta t)dt.
\label{1.7}
\end{equation}
\end{theorem}
The following anysotropic interpolation inequality ``for $m$ points'' was
proved in \cite{b22}.

\begin{theorem}
\label{t1.2}
Let $\overline{p}^j=(p_1^j,\dots,p_n^j)$, $1\leq p_i^j\leq
\infty$, $i=1,\dots,n$, $j=1,\dots,m$, $m\leq n+1$, and $
f\in\bigcap\limits_{j=1}^m L_{\overline{p}^j}$,
$\theta_k\geq 0$, $\sum\limits_{k=1}^m\theta_k=1$. Then $f\in
L_{\overline{p}}$ for any  $\overline{p}$ such
that
$1/{p_i}=\sum\limits_{j=1}^m{\theta_j}/{p_i^j},\quad
i=1,\dots, n$, and
\begin{equation}
\|f\|_{\overline{p}}\leq\|f\|_{\overline{p}^1}^{\theta_1}
\|f\|_{\overline{p}^2}^{\theta_2}
\dots\|f\|^{\theta_m}_{\overline{p}^m}.
\label{1.8}
\end{equation}
\end{theorem}

\subsection{Anysotropic potentials}

Let ${\cal K}_{\overline{\alpha}}(x)$ be a smooth
 function in ${\Bbb R}^n\setminus\{0\}$ such that
${\cal K}_{\overline{\alpha}}(t^\lambda x)=t^{\alpha^*-n}{\cal
K}_{\overline{\alpha}}(x)$, where $x\in {\Bbb R^n}$,
$\lambda\in {\Bbb R^n_+}$, and
%$\lambda_i>0$,
$t>0$ (that is, the function ${\cal K}_{\overline{\alpha}}(x)$
is $\lambda$-homogeneous of order
$\alpha^*-n$).
%We will also consider the kernels, $\lambda$-homogeneous
%in the distribution sense:

%\[
%t^{-n}\langle{\cal K}_{\overline{\alpha}}(x/t^\lambda),\omega\rangle=
%\langle{\cal K}_{\overline{\alpha}},\omega(t^\lambda x)\rangle,\
%\omega\in\Phi_0.
%\]
An operator of the form
\begin{equation}
(K^{\overline{\alpha}}\varphi)(x)=
\int\limits_{{\Bbb R}^n}{\cal K}_{\overline{\alpha}}(x-y)\varphi(y)dy
\label{2.2}
\end{equation}
is called {\it anysotropic potential} (of order $\alpha^*$). This operator
%(\ref{2.2})
is well-defined on the whole space
$L_{\overline{p}}$,
${1\leq p_i\leq\infty}$ ($i=1,\dots,n$),
if $\sum\limits_{i=1}^n\frac{1}{\alpha_ip_i}>1$.
%${1\leq p_i\leq\infty}$ ($i=1,\dots,n$).
The following theorem
can be regarded as an anysotropic analogue of the well-known Sobolev theorem
for the Riesz potential operator.
%provides $L_{\overline{p}}\rightarrow L_{\overline{q}}$
%-estimates for this operator.

\begin{theorem}[\cite{b12}]\label{t2.1}
Let $1< p_i\leq q_i\leq\infty$ ($i=1,\dots, n-1$),
$1<p_n<q_n<\infty$, $\sum\limits_{i=1}^n\frac{1}{\alpha_ip_i}>1$,
and
\begin{equation}
  \sum\limits_{i=1}^n\frac{1}{\alpha_iq_i}=
  \sum\limits_{i=1}^n\frac{1}{\alpha_ip_i}-1.
\label{2.3}
\end{equation}
Then the operator   $K^{\overline{\alpha}}$ is bounded from
$L_{\overline{p}}$ into $L_{\overline{q}}$.
\end{theorem}

\begin{remark}\label{r3.1}
The statement of Theorem~\ref{t2.1} can be generalized to the case of
%If the kernel  ${\cal K}_{\overline{\alpha}}$ in  (\ref {2.2})
  $\lambda$-homogeneous kernels of some order $\gamma<0$
(not necessarily $\gamma=\alpha^*-n$).
One can obtain such a generalization replacing (\ref{2.3})
by the following condition
%then to provide the boundedness of $K^{\overline{\alpha}}$
%we have to replace
%(\ref{2.3}) by the more general condition

\[
  \sum\limits_{i=1}^n\frac{1}{\alpha_iq_i}=
  \sum\limits_{i=1}^n\frac{1}{\alpha_ip_i}-\frac{\gamma+n}{\alpha^*}.
\]

\end{remark}

For $\sum\limits_{i=1}^n\frac{1}{\alpha_ip_i}\leq 1$ we treat the potential
$K^{\overline{\alpha}}\varphi$, $\varphi\in L_{\overline{p}}$, in the sense of $\Phi_0'$-distributions:

\begin{equation}
 \langle K^{\overline\alpha}\varphi,\omega\rangle=
 \langle\varphi,K^{\overline\alpha}\omega\rangle,\quad \omega\in \Phi_0,\quad
 \varphi\in L_{\overline{p}}.
\end{equation}
Such a treatment is correct since
${\cal K}_{\overline{\alpha}}(x)$ is a convolute in $\Phi_0$, in accordance
with the Gel'fand-Shilov theorem (see \cite{Gel}, P.155),
because $\widehat{{\cal K}_{\overline{\alpha}}}(\xi)$ is a multiplier in $\Psi_0$.

We  consider some  anysotropic potentials of special-type.
Let
$Q^{\overline{\alpha}}$ be the operator
with the kernel $Q_{\overline{\alpha}}(x)= F^{-1}
\left(\frac{1}{\sa(\xi)}\right)(x)$, where  $\sa(\xi)$ is a symbol of HSI
(\ref{new:1.3}). The kernel $Q_{\overline{\alpha}}(x)$ admits the following
representations (see \cite{b4} and \cite{b7} in the  cases $\alpha^*<n$ and $\alpha^*\geq n$ respectively):

\begin{eqnarray}
     Q_{\bar\alpha}(x)&=&\int\limits_0^\infty t^{n-\alpha^*-1}
  \tilde g(t^\lambda x)dt\quad , \mbox{ пif }  \quad  \alpha^*<n,
  \label{2.4}\\
    Q_{\bar\alpha}(x)&=&m\sum_{k=|m|}\frac{x^k}{k!}
  \int\limits_0^1   t^{n-\alpha^*+\lambda\cdot k-1}dt
  \int\limits_0^1  (1-u)^{m-1}(D^k\tilde g)(ut^\lambda x)du+
  \nonumber\\
 && + \int\limits_1^\infty  t^{n-\alpha^*-1}
  \tilde g(t^\lambda x)dt,\quad m=[\gamma]+1,   \nonumber\\
 &&\quad  \gamma=(\alpha^*-n)
  \max_j\frac{\alpha_j}{\alpha^*}\quad , \mbox{ if } \quad \alpha^*\geq
  n,   \nonumber
\end{eqnarray}
where $g(\xi)=\varphi(\rho(\xi))/S_{\overline{\alpha}}(\xi)$,
$\varphi\in C_0^\infty(0,\infty)$, $\varphi(t)\geq 0$, and
$\displaystyle \int\limits_0^\infty\frac{\varphi(t)}{t}dt=1$.
%the inverse Fourier transform  is interpreted in the $\Phi'_0$-
%sense.
The function $Q_{\overline{\alpha}}(x)$ is $\lambda$-homogeneous
of order $\alpha^*-n$, if $\alpha ^*<n$.
In the case $\alpha^*\geq n$ it is $\lambda$-homogeneous
in the $\Phi _0'$-sense, that is,
%if $\alpha ^*\geq n$.
\[
t^{\alpha^*-n}\langle\ Q_{\overline{\alpha}}
(x/t^\lambda),\omega\rangle=
\langle\ Q_{\overline{\alpha}},\omega(t^\lambda x)\rangle,\quad \omega\in\Phi_0.
\]

We also consider the modified
potential $\cal R^{\overline{\alpha}}\varphi$ with the kernel
${\cal R}_{\overline{\alpha}}(x)$, being $\lambda$-homogeneous (in the regular sense) for any
$\alpha^*>0$, which differs from $Q_{\overline{\alpha}}(x)$
% $\alpha ^*\geq n$
by a polynomial in the case $\alpha ^*\geq n$ (see
\cite{b7}):
\begin{equation}
\begin{array}{lcc}
{\cal R}_{\overline{\alpha}}(x)=Q_{\overline{\alpha}}(x) ,&
   \mbox{ при } & \mbox{if} \quad \alpha ^*<n,\\
 {\cal R}_{\overline{\alpha}}(x)=
   Q_{\overline{\alpha}}(x)+\sum\limits_{|k|\leq m-1}
\frac{x^k}{k!}\,\frac{\left( D^k\widetilde{g}\right)(0)}
{\lambda\cdot k+n-\alpha ^*},& \mbox{ при } & \mbox{if} \quad  \alpha ^*\geq n.
\end{array}
\label{2.5}
\end{equation}

We observe that
\begin{equation}
{\cal R}^{\overline{\alpha}}\varphi=Q^{\overline{\alpha}}\varphi,\quad \varphi\in \lp ,
\end{equation}
in the  case $\alpha ^*\geq n$, when both potentials are treated in the $\Phi _0'$-sense,
in view of  the fact that $\Phi _0$-functions are orthogonal to all polynomials.

\setcounter{equation}{0}
\setcounter{theorem}{0}
\setcounter{lemma}{0}

\section{Some auxiliary statements}\label{sect2}

\subsection{Characterization of the space $\lpr$ in terms of anysotropic potentials}

Let $Q^{\overline{\alpha}}(\lp)$ be the space of anysotropic potentials:
\[
Q^{\overline{\alpha}}(L_{\overline{p}})=\left\{ f:
f=Q^{\overline{\alpha}}\varphi, \varphi \in
\lp \right\},
\]
\[
\overline{\alpha}=(\alpha_1,\dots,\alpha_n),\quad\overline{p}=(p_1,\dots,p_n),\quad
\alpha _i>0,\quad 1\leq p_i<\infty;
\]
we set
\[
\|f\|_{Q^{\overline{\alpha}}(L_{\overline{p}})}=\|\varphi\|_{{\overline{p}}}.
\]

\begin{theorem}\label{t3.1}
Let $1<p_i<\infty$, $1\leq r_i<\infty$, $\alpha_i>0$, $i=1,\dots
,n$. Then
\[
  \lpr=\lr\cap Q^{\overline{\alpha}}(L_{\overline{p}}).
\]
Moreover,the norms $\displaystyle \|f\|_{\lpr}$ and
$\|f\|_{\overline{r}}+\|f\|_{Q^{\overline{\alpha}}(L_{\overline{p}})}$
are equivalent.
\end{theorem}

\noindent{\bf Proof}. Let $f\in L_{\bar p,\bar
r}^{\bar\alpha}$. We  observe that
\begin{equation}
f=Q^{\bar\alpha}B\varphi,
\label{3.1}
\end{equation}
where
$\varphi=F^{-1}\left(\sum\limits_{j=1}^n|\xi_j|^{\alpha_j}\right)Ff\in
L_{\bar p}$, $B$ is the bounded operator  in $L_{\bar
p}$, generated by the $\overline p$-multiplier
$b(\xi)=S_{\bar\alpha}(\xi)/\sum\limits_{j=1}^n|\xi_j|^{\alpha_j}$
(the relation $b(\xi)\in M_{\overline p}$ is verified with the aid of
 Lizorkin theorem , see \cite{b12}).  Indeed, for $\omega\in\Phi$ we have
\[
\langle f,\omega\rangle=\frac{1}{(2\pi)^n}\langle\hat f,\hat\omega\rangle=
\frac{1}{(2\pi)^n}\langle\hat \varphi,
b(\xi)\frac{\hat\omega(\xi)}{S_{\bar\alpha}(\xi)}\rangle=
\langle B\varphi,Q^{\bar\alpha}\omega\rangle.
\]
From here we derive equality (\ref{3.1}), which yields $f\in Q^{\overline{\alpha}}(\lp)$.
Relatively, $L_{\bar p,\bar
r}^{\bar\alpha}\subset L_{\bar r}\cap
Q^{\bar\alpha}(L_{\bar p})$ and
$\|f\|_{\overline{r}}+\|f\|_{Q^{\overline{\alpha}}(L_{\overline{p}})}\leq
C\displaystyle \|f\|_{\lpr}$.

Let now $f\in L_{\bar r}$ and
$f=Q^{\bar\alpha}\varphi$, $\varphi\in L_{\bar p}$. Then\\
$F^{-1}\left(\sum\limits_{j=1}^n|\xi_j|^{\alpha_j}\right)Ff=B^{-1}\varphi\in
 L_{\bar p}$. Therefore $f\in L_{\bar p,\bar r}^{\bar\alpha}$ and
$\displaystyle \|f\|_{\lpr}\leq C(\|f\|_{\overline{r}}+\|f\|_{Q^{\overline{\alpha}}(L_{\overline{p}})}).
\hfill\blacksquare$

\begin{corol}\label{c4.2}
The spaces $\lpr$ are complete.
\end{corol}

\subsection{On a continuous imbedding of the spaces $\lpr$ with respect to $\overline{\alpha}$}

\begin{theorem}\label{t3.3}
Let  $1\leq r_i<\infty$, $1<p_i<\infty$, $\alpha_i>0$, and
$\beta_i>0$ be such that   $0<\alpha^*-\beta^*<n$ and
 $\displaystyle\sum\limits_{i=1}^n\frac{1}{\alpha_i
p_i}>1-\frac{\beta^*}{\alpha^*}$. Then
\begin{equation}
\lpr\longrightarrow L_{\overline{s},\overline{r}}^{\overline{\beta }},
\label{3.7}
\end{equation}
 $\overline{s}=(s_1,\dots ,s_n)$ being an arbitrary vector such that  $s_i\geq p_i$, if
$i=1,\dots ,n-1$, $s_n>p_n$, and
\[
 \sum\limits_{i=1}^n\frac{1}{\alpha_i s_i}=
\sum\limits_{i=1}^n\frac{1}{\alpha_i
p_i}-\left(1-\frac{\beta^*}{\alpha^*}\right).
\]
\end{theorem}

\noindent{\bf Proof}. Let us consider the anysotropic potential
$I^{\bar\alpha}\varphi=k_{\bar\alpha}*\varphi$ with kernel
$k_{\bar\alpha}(x)\stackrel{(\Phi'_0)}{=}F^{-1}(\rho^{-\alpha^*}(\xi))(x)$,
which is possessed of the same properties, as the kernel $Q_{\overline\alpha}(x)$ (see (\ref{2.4})).

%The convolution is interpreted in the usual sense (i.e. as an
%integral on ${\Bbb R}^n$) in the case
%$\dst\sum_{i=1}^n\frac{1}{\alpha_i p_i}>1$ and in the sense of the
%space $\Phi'_0$ in the case $\dst
%1-\frac{\beta^*}{\alpha^*}<\sum_{i=1}^n\frac{1}{\alpha_i p_i}\leq
%1$ (the function $k_{\bar\alpha}(x)$ is $\lambda$-homogeneous of
%the order $\alpha^*-n$ as $\alpha^*<n$, it is
%$\lambda$-homogeneous in the sense of the space $\Phi'_0$ as
%$\alpha*\geq n$ and $k_{\alpha}(x')\ne0$ in view of the Lemma on
%Fourier transform of $\lambda$-homogeneous function, see
%\cite{b4}). Then

Since both functions $\rho^{\alpha^*}(\xi)/S_{\overline\alpha}(\xi)$ and
$\rho^{-\alpha^*}(\xi)S_{\overline\alpha}(\xi)$ belong to
 $M_{\overline{p}}$
by the Lizorkin theorem (see \cite{b12}), we have

$$
Q^{\bar\alpha}(L_{\bar p})=I^{\bar\alpha}(L_{\bar p})
\stackrel{\rm def}{=}\left\{
f:\  f=I^{\bar\alpha}\varphi,\ \varphi\in L_{\bar p}\right\}.
$$
% and $\|f\|_{Q^{\overline{\alpha}}(L_{\overline{p}})\eqiv\|f\|_{I^{\overline{\alpha}}(L_{\overline{p}})$.
Therefore
\begin{equation}
  L_{\bar p,\bar r}^{\bar\alpha}= L_{\bar r}\cap I^{\bar\alpha}(L_{\bar p})
\label{new:4.3}
\end{equation}
 by
Theorem~\ref{t3.1}.

Let $f=I^{\bar\alpha}\varphi$, $\varphi\in L_{\bar p}$, and $\omega\in \Phi_0$. We have
\[
\langle f,\omega\rangle=\frac{1}{(2\pi)^n}\langle
\hat\varphi,\rho^{-\alpha^*}(\xi)\hat\omega(\xi)
\rangle=\langle\psi,I^\beta\omega\rangle,
\]
where $\psi=F^{-1}(\rho^{\beta^*-\alpha^*}(\xi))*\varphi$.
Making use of Remark \ref{r3.1}, we easily obtain the imbedding
$Q^{\overline{\alpha}}(L_{\overline{p}})\rightarrow I^{\overline{\beta}}(L_{\overline{s}})$,
which yields (\ref{3.7}) by (\ref{new:4.3}).

%obtain $\|\psi\|_{\bar s}\leq c\|\varphi\|_{\bar p}$ by the
%theorem on a boundedness of potential-type operators with
%homogeneous kernels (see \cite{b12}). Thus
%$\|f\|_{I^{\bar\beta}(L_{\bar s})}\leq
%c\|f\|_{I^{\bar\alpha}(L_{\bar p})}$ which gives the embedding
%(\ref{3.6}).\quad
$\hfill\blacksquare$

\setcounter{equation}{0}
\setcounter{theorem}{0}
\setcounter{lemma}{0}
\section{Proof of the main results}\label{new:sect5}

\subsection{Proof of Theorem~\ref{new:t3.1}}

Let $f\in\lpr$, then
$f=Q^{\bar\alpha}\varphi$, $\varphi\in\lp$,
by  Theorem~\ref{t3.1}. Relatively, $f={\cal
R}^{\bar\alpha}\varphi$. We prove the following integral representation
\begin{equation}
\left(T_\varepsilon^{\bar\alpha}{\cal R}^{\bar\alpha}\varphi\right)(x)=
\int\limits_{{\Bbb R}^n}{\cal K}(y)\varphi(x-\varepsilon^\lambda y)dy,
\label{3.3}
\end{equation}
where the kernel
\begin{equation}
 {\cal K}(y)=\rho^{-n}(y) \int\limits_{\rho(t)>1/\rho(y)}
 \frac{(\Delta_t^{2\ell}{\cal R}_{\bar\alpha})
 \left(\frac{y}{\rho^\lambda(y)}\right)}{\rho^{n+\alpha^*}(t)}dt
\label{3.4}
\end{equation}
is  averaging. We note that the representation (\ref{3.3}) was proved in \cite{b7} for $\varphi\in\Phi_0$.

Since the symbol
\[
 S_{\bar\alpha,\varepsilon}(\xi)= (-1)^\ell
2^{2\ell}\!\!\int\limits_{\rho(t)>\varepsilon}\!\!
\frac{\sin^{2\ell}\frac{\xi\cdot t}{2}}{\rho^{n+\alpha^*}(t)}dt
\]
of the "truncated" integral (\ref{new:1.4})
is a multiplier in $\Psi_0$ (see \cite{b11}), we conclude
that the space $\Phi_0$ is invariant with respect to the operator $T_{\varepsilon}^{\bar\alpha}$.
Therefore
\begin{eqnarray*}
 \langle T_{\varepsilon}^{\bar\alpha}f,\omega\rangle&=&
 \langle f,T_{\varepsilon}^{\bar\alpha}\omega\rangle=
 \langle\varphi,{\cal R}^{\overline{\alpha}}T_{\varepsilon}^{\bar\alpha}\omega\rangle=
 \langle\varphi,T_{\varepsilon}^{\bar\alpha}{\cal R}^{\bar\alpha}\omega\rangle=
 \langle \varphi,\int\limits_{{\Bbb R}^n}{\cal K}(y)\omega(x-\varepsilon^\lambda y)dy
\rangle= \\
&=&\langle\int\limits_{{\Bbb R}^n}{\cal
K}(y)\varphi(x-\varepsilon^\lambda y)\,dy,\omega\rangle,\quad \omega\in\Phi_0.
\end{eqnarray*}
Thus (\ref{3.3}) is valid in the $\Phi'_0$-sense.

Since two locally integrable functions, which agree in the $\Phi'_0$-sense,
may differ from each other only by a polynomial, we obtain (\ref{3.3}) for almost all
$x\in{\Bbb R}^n$ ( because both sides of (\ref{3.3}) belong to the weighted spaces
$L_{1,\gamma}=\lbrace f(x):\dst\int\limits_{{\Bbb R}^n}{\mid f(x)\mid (1+\mid x\mid)}^\gamma dx<\infty\rbrace$
for some $\gamma>-n$).
% We observe, that two local summable functions which coincides in the
%sense of the space $\Phi'_0$ may differ themselves only by a
%polynomial. Since $T_{\varepsilon}^{\bar\alpha}f\in L_{\bar
%r}({\Bbb R}^n)$ and $\dst\int\limits_{{\Bbb R}^n}{\cal
%K}(y)\varphi(x-\varepsilon^\lambda y)dy\in L_{\bar p}({\Bbb R}^n)$
%does not contain a polynomial, the representation (\ref{3.3}) is
%true for almost all $x\in {\Bbb R}^n$.

Letting $\varepsilon\to0$ in (\ref{3.3}), in view of  Theorem \ref{t1.1} we have
\[
T^{\bar\alpha}f=\varphi\in \lp,\quad
\|T^{\bar\alpha}f\|_{\bar p}=\|f\|_{Q^{\bar\alpha}(L_{\bar p})}.
\]

Let now $f\in\lr$, $T^{\bar\alpha}f\in\lp$. The statement of  theorem will follow from the relation
$f=Q^{\bar\alpha}\varphi$, where $\varphi=T^{\bar\alpha}f$. For
$\omega\in\Phi_0$ we have
\begin{eqnarray}
 \langle \varphi,Q^{\bar\alpha}\omega\rangle &=&
 %\langle\varphi,{\cal R}^{\bar\alpha}\omega\rangle =
\langle \lim_{\varepsilon\to0\atop{(L_{\bar p})}}T_\varepsilon^{\bar\alpha}f,
{\cal R}^{\bar\alpha}\omega\rangle =
      \lim_{\varepsilon\to0}\langle T_\varepsilon^{\bar\alpha}f,
      {\cal R}^{\bar\alpha}\omega\rangle =
      \nonumber\\
% &=&\lim_{\varepsilon\to0}\langle f,T_\varepsilon^{\bar\alpha}
% {\cal R}^{\bar\alpha}\omega\rangle =\lim_{\varepsilon\to0}\langle T_\varepsilon^{\bar\alpha}f,
%      {\cal R}^{\bar\alpha}\omega\rangle = \nonumber\\
&=&       \lim_{\varepsilon\to0}\langle f,\int\limits_{{\Bbb R}^n}
 {\cal K}(y)\omega(x-\varepsilon^\lambda y)dy\rangle =\langle
f,\omega\rangle .
\label{3.5}
\end{eqnarray}
The second of equalities of this chain follows from the fact
that the convergence in $L_{\bar p}$  implies that in
$\Phi_0'$, the last one is justified by the H\"older
inequality. Thus we have proved that
\begin{equation}
  \langle f,\omega\rangle=\langle \varphi,Q^{\overline{\alpha}}
  \omega\rangle.
\label{3.6}
\end{equation}
In the  case $\dst\sum_{i=1}^n\frac{1}{\alpha_ip_i}> 1$ we have
\[
 \langle f,\omega\rangle=\langle Q^{\overline{\alpha}}\varphi,
  \omega\rangle
\]
hence, $f=Q^{\bar\alpha}\varphi$ for almost all $x\in {\Bbb R}^n$. In  the case
$\dst\sum_{i=1}^n\frac{1}{\alpha_ip_i}\leq 1$ the relation (\ref{3.6}) itself means
that $f=Q^{\bar\alpha}\varphi$. Now the application of Theorem~\ref{t3.1}
yields $f\in \lpr$.\quad
$\hfill\blacksquare$

%\subsection{On the continuity modulus of functions in $\lpr$}

%We need  the following properties of  finite differences of  functions in $f\in\lpr$.

\begin{corol}
%[finite differences' representation]
\label{l5.1}
Let $f\in\lpr$,
$\alpha_i>0$, $1<p_i<\infty$, $1\leq r_i<\infty$, $i=1,\dots,n$.
%and let the restrictions (\ref{2.6}) are valid.
Then
%for $\dst m>\max_j\alpha_j$, $t\in{\Bbb R}^n$
the following integral representation holds:
\begin{equation}
(\Delta_t^m f)(x)=\int\limits_{{\Bbb R}^n}(\Delta_t^m
 {\cal R}_{\bar\alpha})(x-y)(T^{\bar\alpha}f)(y)dy,\quad \dst m>\max_j\alpha_j.
\label{5.21}
\end{equation}
\end{corol}

\begin{corol}\label{l5.2}
Let $f\in\lpr$,
$\alpha_i>0$, $1<p_i<\infty$, $1\leq r_i<\infty$ ($i=1,\dots,n$).
%s (\ref{2.6}) are valid.
Then
\begin{equation}
\|\Delta_t^m f\|_{\bar p}\leq c\rho^{\alpha^*}(t)
\|T^{\bar\alpha}f\|_{\bar p},\quad \dst m>\max_j\alpha_j.
\label{5.22}
\end{equation}
\end{corol}
We give the proofs of Corollaries~\ref{l5.1} and~\ref{l5.2} in
Appendix.

\subsection{Proof of Theorem~\ref{new:t4.1}}
We first  approximate the function $f\in\lpr$ by smooth functions
using the averages (\ref{1.7}) with a "very nice" kernel.
Namely, we assume that $a(t)\in C_0^\infty$. Then $\displaystyle
\|f-f_\delta\|_{\lpr}\to0$, as $\delta \to0$. Therefore it remains to
approximate $f\in C^\infty \cap\lpr$ by $C_0^\infty$-functions.

For this goal, we consider the function $\mu (x)\in C_0^\infty$  supported in the
ball $\rho (x)<2$ such that $0\leq\mu(x)\leq1$ , $\mu(x)=1$ if $\rho
(x)<1$, and $\mu(x)=0$ in the case $\rho(x)\geq2$.
%$|\mu(x)|\leq1$.
Let $\mu _N(x)=\mu
(N^{-\lambda }x)$. We  show that the smooth truncations
$f_N(x)=\mu _N(x)f(x)(\in C_0^\infty)$ approximate $f(x)$ in the norm of $\lpr$ , as
$N\to\infty $. Let $\nu (x)=1-\mu (x)$, $\nu
_N(x)=\nu(x/N)$. We have to prove that
\begin{equation}
\left\|\ta(\nu_Nf)\right\|_{\overline{p}}\to0,\quad as\quad N\to\infty,\quad
f\in C^\infty\cap\lpr ,
\label{4.1}
\end{equation}
where
\begin{equation}
  (\ta f)(x)=
 \int\limits_{{\Bbb R}^n}\frac{(\Delta_t^{2\ell}f)(x)}
{\rho ^{n+\alpha ^*}(t)}\,dt,\quad 2\ell>\max_j\alpha_j.
\label{2.1}
\end{equation}
It should be noted that the integral (\ref{2.1}) converges absolutely
for the functions $f(x)$  bounded in ${\Bbb R}^n$ together with all their  partial derivatives up to the
order $[\alpha^*]+1$.

%introduced by P.\,I.~Lizorkin \cite{b11}.
% We will
%interpret the hypersingular integral $T^{\overline{\alpha}}$ as
%convectionally convergent at origin in $\lp$, for the functions
%$f(x)$, having representation in the form of anysotropic
%potentials with the density from $\lp$:

{\bf I. Proof of (\ref{4.1}) in the case $\displaystyle
\sum\limits_{i=1}^{n}\frac{1}{\alpha _ip_i}>1$. }

We first  suppose that $0 < \max\limits_j \alpha_j < 1$ and
take $\ell=1$ in (\ref{2.1}). Then
\begin{equation}
    \left|\ta\left(\nu_N f\right)(x)\right| \leq
\left| \nu_N(x) \left(\ta f\right)(x) \right| + 2
\left|\left(B_N f\right)(x)\right|,     \label{4.2}
\end{equation}
where
\[
\displaystyle \left( B_N f\right)(x) = \int\limits_{{\Bbb
R}^n}\frac{\nu_N(x-t) - \nu_N(x)}{\rho^{n+\alpha^*}(t)} f(x-t)\,
dt.
\]
 In order to obtain (\ref{4.1}), it  suffices to verify
that
\begin{equation}
    \left\|B_N f\right\|_{\overline{p}} \to 0, \qquad N\to \infty
\label{4.3}
\end{equation}
(since $\|\nu_N(x)f(x)\|_{\overline{p}}\rightarrow 0$, as $N\rightarrow\infty$,
for every $f(x)\in\lp$).

We justify (\ref{4.3})  for $f\in
L_{\overline{q}}$, where $\overline{q}$ is an arbitrary vector
such that ${\overline{q}}>{\overline{p}}$ and equality (\ref{2.3}) is fulfilled.
%$q_i>p_i$ $(i=1,\dots,n)$. Then
This implies (\ref{4.3}) in view of the imbedding
%for $f\in\lpr$ will follow from the imbedding
$\lpr\subset L_{\overline{q}}$, which is
valid by virtue of Theorems~\ref{t3.1} and \ref{t2.1}.

Owing to the uniform estimate
\begin{equation}
\left\|B_N f\right\|_{\overline{p}} \leq C\|f\|_{\overline{q}}, \label{4.4}
\end{equation}
with the constant $C$ not depending on $N$ (see Appendix), in accordance with the
Banach--Steinhaus theorem, it remains to verify (\ref{4.3})
% $L_{\overline{q}}$ set $C_0^{\infty}$.
for $f \in C_0^{\infty}$ (since the class $C_0^{\infty}$ is dense in
$L_{\overline{q}}$).

Assuming that $f(x) \in C_0^{\infty}$  is supported  in the
ball $\rho(x)<a$, we have
\begin{eqnarray*}
&& \left\| B_N f\right\|_{\overline{p}} \leq
   \frac{C}{N^{\theta}} +
   CN^{\alpha^{\ast}\sum\limits_{i=1}^n
   \frac{1}{\alpha_i p_i}-\alpha^{\ast}}
   \int\limits_{\rho(t) > \frac{2a}{C_1N}}
   \frac{dt}{\rho^{n+\alpha^{\ast} - \theta}(t)} \times \\
&&  \qquad \times \left( \int\limits_{{\Bbb R}^1}
dx_n \dots \left( \int\limits_{{\Bbb R}^1} \frac{ |f(N^{\lambda}
x)|^{p_1} dx_1}{(1+\rho(x))^{p_1 \theta}(1+
\rho(x-t))^{p_1\theta}}
\right)^{\frac{p_2}{p_1}} \dots \right)^{\frac{1}{p_n}},
\end{eqnarray*}
where $\theta=\alpha^*/\max\limits_{j}\alpha_j$, $C_1$ being such that $\rho(x-t) \geq C_1 \rho(t) - \rho(x)$.
Since $\displaystyle \rho(x-t) \geq C_1 \rho(t) - \frac{a}{N}
>0$, if $\displaystyle \rho(x) < \frac{a}{N}$ and $\displaystyle
\rho(t) > \frac{2a}{C_1N}$, we have
\[
    \left\|B_N f\right\|_{\overline{p}} \leq \frac{C}{N^{\theta}} +
\frac{C}{N^{ \alpha^{\ast} }} \int\limits_{{\Bbb R}^n}
\frac{dt} {\rho^{n+\alpha^{\ast}-\theta}(t)\left( 1+ C_1\rho(t) -
\frac{a}{N_0} \right)^{\theta}} \to 0,
\]
as $N\to \infty$ $(N>N_0>a)$. Thus we have obtained the desired result
in the case $0< \max\limits_j \alpha_j <1$.

In order to prove (\ref{4.1}) in the  case $\max\limits_j \alpha_j \geq 1$,
we use the  induction argument. It should be noted that such an idea
%to prove the denseness of the class  $C_0^{\infty}$ in $ L_{p,r}^\alpha$
is due to S.G. Samko (see \cite{b17}) in the isotropic case.
%is due to S.G. Samko (see \cite{b17}).
Its realization in the most general anysotropic case of vector-valued
$\overline{\alpha}$, $\overline{p}$,
and $\overline{r}$, treated here, is much more difficult. The following lemma,
which enables one us to apply the induction argument, plays a
crucial role in the proof of Theorem \ref{new:t4.1}.

\begin{lemma}\label{l4.1}
Let $m=1,2,\dots$ and $\displaystyle \sum\limits_{i=1}^n \frac{1}{\alpha_i p_i} >1$.
If the class $C_0^{\infty}$ is dense in $\lpr$
for $0<\max\limits_j \alpha_j < m$,
% and vectors $\overline{p}$,
%$\overline{r}$ and $\overline{\alpha}$ with components
%$1<p_i<\infty$, $1\leq r_i<\infty$, $\alpha_i>0$ $(i=1,\dots,n)$,
then it is also dense in $\lpr$ for $m \leq \max\limits_j \alpha_j
< m+1$.
\end{lemma}

\noindent{\bf Proof.} As above, we have to check (\ref{4.1}) for
$m\leq\max\limits_j \alpha_j < m+1$.

Owing to the evident formulas $\displaystyle \left(
\Delta_t^{2\ell}
\nu_N f
\right)(x) = \left( \tilde{\Delta}_t^{2\ell} \nu_N f
\right)(x+\ell t)$,\\
 $\displaystyle \left( \tilde{\Delta}_t^{2\ell}
\nu_N f \right)(x) = \sum\limits_{k=0}^{2\ell} C_{2\ell}^k \left(
\tilde{\Delta}_t^k \nu_N \right)(x) \left(
\tilde{\Delta}_t^{2\ell-k} f \right)(x-kt)$, we have
\[
  T^{\overline{\alpha }} (\nu_N f) = B_N f +
\sum\limits_{k=1}^{2\ell} C_{2\ell}^k B_{N,k} f,
\]
where
\begin{eqnarray*}
    &&\left(B_N f\right)(x) = \int\limits_{{\Bbb R}^n} \frac{\nu_N(x+\ell t)
\left( \tilde{\Delta}_t^{2\ell} f \right)(x+\ell
t)}{\rho^{n+\alpha^{\ast}}(t)}\, dt; \\
    &&\left( B_{N,k} f\right)(x) = \int\limits_{{\Bbb R}^n} \frac{\left(
\tilde{\Delta}_t^k \nu_N \right)(x+\ell t) \left(
\tilde{\Delta}_t^{2\ell-k} f
\right)(x+(\ell-k)t)}{\rho^{n+\alpha^{\ast}}(t)}\, dt,  \\
    && \qquad k = 1,2,\dots, 2\ell.
\end{eqnarray*}
Let us show that
\begin{eqnarray}
&& \left\|B_{N,k} f\right\|_{\overline{p}} \to 0, \quad as \quad  N\to
\infty, \quad k = 1,2, \dots, 2\ell,   \label{4.5}
\end{eqnarray}
and
\begin{eqnarray}
  && \left\|B_N f \right\|_{\overline{p}} \to 0, \quad as \quad N \to
\infty.     \label{4.6}
\end{eqnarray}
In the case $k=2\ell$ the relation (\ref{4.5}) is proved for $f\in
L_{\overline{q}}$  $(\supset L_{\overline{p},
\overline{r}}^{\overline{\alpha }})$ in just the same way, as (\ref{4.1})
for $0 < \max\limits_j \alpha_j < 1$ .
% where
%$\overline{q}$ is vector with components satisfying (\ref{2.3})
%with $q_i>p_i$. The uniform bound $\displaystyle
%\bigl\|B_{N,2\ell} f
%\bigr\|_{\overline{p}} \leq c \|f\|_{\overline{q}}$, $f \in
%L_{\overline{q}}({\Bbb R}^n)$, is proved in the same way as
%(\ref{4.4}). Moreover, assuming that $f(x) \in C_0^{\infty}$ and
%its support lies in the ball $\rho(x) <a$, and applying the
%procedure of proof (\ref{4.3}), we get $\bigl\|B_{N,2l}
%f\bigr\|_{\overline{p}} \to 0$ as $N \to \infty$, $f \in
%C_0^{\infty}$. Since $C_0^{\infty}$ is dense in
%$L_{\overline{q}}({\Bbb R}^n)$ according to the Banach--Steinhaus
%theorem the relation (\ref{4.5}) is true for $k = 2\ell$.

In the cases $k=1,2,\dots,2\ell-1$ we  prove the
stronger relation:
\begin{equation}
    \left\|B_{N,k} f\right\|_{\overline{p}} \to 0, \quad as \quad N \to
\infty, \quad f \in L_{\overline{\tau }^k,
\overline{r}}^{\overline{\gamma }^k},  \label{4.7}
\end{equation}
 the vectors $\overline{\gamma }^k =
(\gamma_1^k, \dots, \gamma_n^k)$ and $\overline{\tau }^k =
(\tau_1^k, \dots, \tau_n^k)$ being such that  $\displaystyle
\sum\limits_{i=1}^n \frac{1}{\alpha_i
\tau_i^k}=\sum\limits_{i=1}^n\frac{1}{\alpha_i p_i} - \left( 1 -
\frac{(\gamma^k)^{\ast}}{\alpha^{\ast}} \right)$ , $\displaystyle
\tau_i^k > p_i$, and
\begin{equation}
0<\max\limits_i \gamma_i^k <m.
\label{new:0}
\end{equation}
(Later we will put some additional restrictions on $\gamma_i^k$).
% $\displaystyle\frac{1}{(\gamma^k)^{\ast}} =
%\frac{1}{n} \sum\limits_{i=1}^n \frac{1}{\gamma_i^k}$,
Then (\ref{4.7}) yields (\ref{4.5}) by (\ref{3.7})
%the application of Theorem~\ref{t3.3} yields (\ref{4.5})
for $f \in L_{\overline{p},
\overline{r}}^{\overline{\alpha }}$ and $\max\limits_i
\gamma_i^k < m \leq \max\limits_i \alpha_i$. We first obtain the
uniform estimate
\begin{equation}
    \left\|B_{N,k} f\right\|_{\overline{p}} \leq C
 \|f\|_{L_{\overline{\tau }^k,
 \overline{r}}^{\overline{\gamma }^k}}, \qquad  f \in
 L_{\overline{\tau }^k,\overline{r}}^{\overline{\gamma }^k}.
\label{4.8}
\end{equation}
Making use of (\ref{1.5}) and the H\"older inequality , we obtain
\begin{eqnarray*}
\bigl\|B_{N,k} f \bigr\|_{\overline{p}}& \leq &
 \frac{C}{N^{k\theta}}\! \int\limits_{{\Bbb R}^n}
 \frac{\left\|\tilde{\Delta}_t^{2\ell-k}\!
 f\right\|_{\overline{\tau }^k}}{\rho^{n+\alpha^{\ast}-k\theta}(t)}
 dt\times \\
 &&\times
 \left(\int\limits_{{\Bbb R}^1} dx_n \dots\left(\,
\int\limits_{{\Bbb R}^1}\frac{dx_1}{\prod\limits_{i=0}^k\left(1+
\frac{\rho(x+(\ell-i)t)}{N}\right)^{\theta\xi_1^k}}
\right)^{\xi_2^k/\xi_1^k}\dots \right)^{1/\xi_n^k}\!\!\!\!\!,
\end{eqnarray*}
where $\displaystyle\xi_i^k = \frac{p_i \tau_i^k}{\tau_i^k - p_i}$.
 Applying Lemma~\ref{l1.3}, under the additional
condition $2l-k > \max\limits_j \gamma_j^k$,
we have
\[
  \left\|\Delta_t^{2l -k} f \right\|_{\overline{\tau }^k} \leq c
\rho^{(\gamma^k)^{\ast}}(t) \left\|T^{\overline{\gamma }^k}
f\right\|_{\overline{\tau }^k}, \quad f \in L_{\overline{\tau }^k,
\overline{r}}^{\overline{\gamma }^k} ,
%({\Bbb R}^n),
\]
where
\begin{equation}
\sum\limits_{i=1}^n \frac{1+j_i}{\gamma_i^k}\ne1, \quad
|j|=0,1,\dots,
\left[ \max\limits_i \gamma_i^k
   \left(1-\sum\limits_{i=1}^n \frac{1}{\gamma_i^k} \right)
\right]-1.
\label{4.9}
\end{equation}
Besides this, we assume that
%taking the coordinates of vector
% $\overline{\gamma }^k$ such that
\begin{equation}
\displaystyle \gamma_i^k<\frac{(2\ell-k)\alpha_i}{\max\limits_j
\alpha_j},\quad i=1,\dots,n,
\label{new:1}
\end{equation}
then
%and (\ref{4.9}) are valid, we get
\begin{eqnarray*}
    && \left\|B_{N,k} f \right\|_{\overline{p}} \leq
C \left\|T^{\overline{\gamma }^k} f \right\|_{\overline{\tau}^k}
\int\limits_{{\Bbb R}^n}
\frac{dt}{\rho^{n+\alpha^{\ast}-k\theta-(\gamma^k)^{\ast}}} \times \\
    && \qquad\times \left(\, \int\limits_{{\Bbb R}^1} dx_n \dots \left(\,
\int\limits_{{\Bbb R}^1} \frac{dx_1}{\prod\limits_{i=0}^k \left[ 1+ \rho(x+(l
-i)t)\right]^{\theta\xi_1^k}} \right)^{\xi_2^k/\xi_1^k}
\dots \right)^{1/\xi_n^k}.
\end{eqnarray*}
In order to apply Lemma~\ref{l1.0}, we also put the following restrictions on
 $\overline{\gamma }^k$:
\begin{equation}
 \left(
\gamma^k\right)^{\ast} \displaystyle > \alpha^{\ast}-\theta k .
%< \alpha^{\ast} - \theta.
\label{new:2}
\end{equation}
Then
%$\displaystyle \frac{\theta}{\alpha^{\ast}} <\sum\limits_{i=1}^n
%\frac{1}{\alpha_i \xi_i^k} < \frac{ \theta k}{\alpha^{\ast}}$
%and in view of (\ref{1.1}) we have
Lemma~\ref{l1.0} yields
\begin{equation}
    \left\|B_{N,k} f \right\|_{\overline{p}} \leq C
\left\|T^{\overline{\gamma }^k} f \right\|_{\overline{\tau }^k}.
\label{4.10}
\end{equation}
Taking into account (\ref{new:0}), (\ref{new:1}), and  (\ref{new:2}),
we arrive at the following inequalities for the coordinates ${\gamma_i^k}$:
% to take  vector $\overline{\gamma }^k$ in such a way that

\begin{equation}
\max \left( 0, \alpha_i - \frac{ \alpha_i k}{\max\limits_j
\alpha_j} \right) < \gamma_i^k <
 \min \left( \frac{
m\alpha_i}{\max\limits_j \alpha_j},
\frac{(2\ell-k)\alpha_i}{\max\limits_j
\alpha_j} \right).
\label{new:3}
\end{equation}
%\alpha_i - \frac{\alpha_i}{\max\limits_j \alpha_j} ,
 Since
$\displaystyle \max_j \alpha_j<2\ell$ and $\displaystyle
m\leq\max_j \alpha_j<m+1$, this interval is not empty.

Thus we have proved (\ref{4.8}).
%the uniform estimate (\ref{4.8}) follows from (\ref{4.10}) under
%the choice $\gamma_i^k$ for which (\ref{new:3})
%and (\ref{4.9}) are valid.
In accordance with the Banach--Steinhaus theorem, it suffices to verify
(\ref{4.7})
 on a dense set in $L_{\overline{\tau }^k,
\overline{r}}^{\overline{\gamma }^k}$.
%({\Bbb R}^n)$.
Since $0 <
\max\limits_i \gamma_i^k < m$, the class $C_0^{\infty}$ is dense in
$L_{\overline{\tau }^k,\overline{r}}^{\overline{\gamma }^k}$
according
to  the assumption of induction. For $f(x) \in C_0^{\infty}$, in view of
(\ref{1.5}) we have
\begin{eqnarray*}
    && \left\|B_{N,k} f\right\|_{\overline{p}} \leq \frac{C}{
N^{(\gamma^k)^{\ast}-(2\ell -k)\theta}} \int\limits_{{\Bbb R}^n}
\frac{dt}{\rho^{n+\alpha^{\ast}-2\ell\theta}(t)(1+N\rho(t))^{(2\ell
-k)\theta}} \times \\
    && \qquad\times \left( \, \int\limits_{{\Bbb R}^1} dx_n \dots \left( \,
\int\limits_{{\Bbb R}^1} \frac{ dx_1}{\prod\limits_{i=0}^k (1+ \rho(x-it))^{\theta
\xi_1^k}}\right)^{\xi_2^k/\xi_1^k} \dots
\right)^{1/\xi_n^k}.
\end{eqnarray*}
Applying the H\"older inequality to the inner integral, we get
\begin{equation}
 \left\|B_{N,k} f \right\|_{\overline{p}} \leq
  \frac{c}{N^{(\gamma^k)^{\ast}-(2\ell-k)\theta+d}}
\int\limits_{{\Bbb R}^n} \frac{dt}{ \rho^{n+\alpha^{\ast}-2\ell \theta+d}(t)
(1+\rho(t))^{(k+1)\theta - \alpha^{\ast} + (\gamma^k)^{\ast}}},
\label{new:4}
\end{equation}
where $0 \leq d \leq (2\ell -k)\theta$. The integral on the right-hand side
converges, if we take $d \in ( (2\ell-k)\theta-(\gamma^k)^{\ast},
\min\{(2\ell-k)\theta, 2\ell\theta-\alpha^{\ast}\})$
(this interval is not empty in view of (\ref{new:2})).
Then the right-hand side of (\ref{new:4}) tends to zero, as $N\to\infty$.

Let us prove (\ref{4.6}). We have
\[
\displaystyle \left(B_N
f\right)(x)=\nu_N(x) \left( T^{\overline{\alpha }}
f\right)(x)+(M_N f)(x),
\]
 where
\begin{equation}
    \left(M_N f\right)(x) = \int\limits_{{\Bbb R}^n} \frac{\nu_N(x+\ell t) -
\nu_N(x)}{\rho^{n+\alpha^{\ast}}(t)} \, \left(
\tilde{\Delta}_t^{2\ell} f\right) (x+\ell t)\, dt.  \label{4.11}
\end{equation}
 The validity of (\ref{4.6})  follows from the relation
\[
\left\|M_N f \right\|_{\overline{p}} \to 0,\quad as\quad  N \to \infty,
\quad f \in L_{\overline{p},\overline{r}}^{\overline{\alpha }},
\]
which is verified in just the same way, as (\ref{4.5}).

{\bf II. The case $\displaystyle
\sum\limits_{i=1}^n\frac{1}{\alpha_i p_i} \leq 1$.}

In the case $0<\max\limits_j\alpha _j<1$, as above, we have (\ref{4.2}).

Let $p_j<r_j$ for $j=m_1^1,\dots ,m_{k_1}^1$ and $p_j\geq r_j$ for
$j=m_1^2,\dots ,m_{k_2}^2$, and let $\overline{r}^1$ be the vector
with the following coordinates: $r^1_j=r_j$ for $j=m_1^1,\dots ,m_{k_1}^1$
and $r^1_j=p_j$ for $j=m_1^2,\dots ,m_{k_2}^2$. We observe that
$f\in L_{\overline{r}^1}$
%({\Bbb R}^n)$
by virtue of (\ref{1.8}). Similarly to (\ref{4.3}) we prove that
\begin{equation}
\left\| B_Nf\right\|_{\overline{p}}\to0,\quad as\quad N\to\infty,\quad f\in L_{\overline{r}^1}.
\label{4.12}
\end{equation}

Making use of the Minkowsky and H\"older inequalities ,
we arrive at the estimate
\begin{eqnarray*}
 \|B_N f\|_{\overline{p}}&\leq & N^{-\alpha ^*+\alpha
^*\sum\limits_{j=m_1^1}^{m_{k_1}^1}\frac{1}{\alpha _jp_j}-
\alpha ^*\sum\limits_{j=m_1^1}^{m_{k_1}^1}\frac{1}{\alpha _jr_j}}
\|f\|_{\overline{r}^1}
 \int\limits_{{\Bbb R}^n}
\frac{\|\nu(x-t)-\nu(x)\|_{\overline{s}}} {\rho^{n+\alpha
^*}(t)}\,dt,
\end{eqnarray*}
where $\displaystyle s_j=\frac{p_jr_j}{r_j-p_j}$ for
$j=m_1^1,\dots ,m_{k_1}^1$ and $s_j=\infty $ for $j=m_1^2,\dots
,m_{k_2}^2$. Since 
$\displaystyle
\alpha^*-\alpha^*\sum_{j=m_1^1}^{m_{k_1}^1} \frac{1}{\alpha_j
p_j}+\alpha^*\sum_{j=m_1^1}^{m_{k_1}^1} \frac{1}{\alpha_j r_j}>0$
and the integral on the right-hand side converges, we obtain (\ref{4.12}).

In oder to extend the result on denseness of $C_0^\infty$
in $\lpr$ to the case
% Let us extend the feasibility of (\ref{4.12}) for
$\max\limits_{j}\alpha _j\geq1$, we use the following analogue of Lemma~ \ref{l4.1}.

\begin{lemma} \label{l4.2}
Let $m=1,2,\dots $, $1<p_i<\infty $, $1\leq r_i<\infty $, $\alpha
_i>0$, and $\dst \sum_{i=1}^n
\frac{1}{\alpha_i p_i}\leq 1$ . If
$C_0^\infty $ is dense in $\lpr$ for $\displaystyle
0<\max\limits_j\alpha _j<m$, then it is dense in $\lpr$ for
$\displaystyle m\leq\max\limits_j\alpha _j<m+1$.
\end{lemma}

The proof of this lemma is much in lines with that of Lemma~ \ref{l4.1}.
The only difference is as follows. To prove (\ref{4.5}), we have to use
 the assumption of induction only in the  case $\displaystyle \frac{1}{\max\limits_j\alpha _j}
-\sum\limits_{j=m_1^1}^{m_{k_1}^1}\frac{1}{\alpha _jp_j}+
\sum\limits_{j=m_1^1}^{m_{k_1}^1}\frac{1}{\alpha _jr_j}\leq 0$
 and $k=1,\dots,\bigl[\max\limits_j
\alpha_j \bigr]$.
As in the proof of Lemma~\ref{l4.1} , we base ourselves on the imbedding (\ref{3.7}). Namely,
we choose $\beta_i$ $(i=1,\dots,n)$  such  that
 $\alpha^{\ast}\,-\,\beta^{\ast}\,<\,\min\{n,\theta\}$, $0 <
\max\limits_j \beta_j < m$, $\displaystyle \sum\limits_{i=1}^n
\frac{1+j_i}{\beta_i}\ne1$, $\displaystyle|j|=0,1,\dots,
\left[\max\limits_i\beta_i\left( 1-\sum\limits_{i=1}^n
\frac{1}{\beta_i}\right) \right] -1$, and check (\ref{4.5}) for
$f \in L_{\overline{\zeta}, \overline{r}}^{\overline{\beta
}}$, where the components of $\overline{\zeta}$
($\zeta_i>p_i$, $i=1,\dots,n$) satisfy the  equality
$\displaystyle\sum\limits_{i=1}^n\frac{1}{\zeta_i
\alpha_i} =
\sum\limits_{i=1}^n
\frac{1}{\alpha_i p_i} - \left(
1-\frac{\beta^{\ast}}{\alpha^{\ast}}\right)$.

As regards the rest values of $k$, the relation (\ref{4.5}) is proved by the
direct estimation of the norm $\|B_{N,K} f\|_{\overline{p}}$ .
We leave the corresponding details to the reader.\quad $\hfill\blacksquare$

%%%%%%%%%%%%%%%
\setcounter{equation}{0}
\setcounter{theorem}{0}
\setcounter{lemma}{0}

\section{Appendix}\label{sect3}

\subsection{Proof of equality (\ref{5.21}). }

Let $f\in L_{\bar p,\bar
r}^{\bar\alpha}$. We denote
$\varphi_\varepsilon=T_\varepsilon^{\bar\alpha}f$, then
\begin{equation}
\int\limits_{{\Bbb R}^n}(\Delta_t^m
 {\cal R}_{\bar\alpha})(x-y)\varphi_\varepsilon(y)dy=
\int\limits_{{\Bbb R}^n}(\Delta_t^m f)(x-\varepsilon^\lambda y)
 {\cal K}(y)dy,
\label{new:5}
\end{equation}
where ${\cal K}(y)$ is the kernel (\ref{3.4}). Let us pass to the
limit, as $\varepsilon\to0$, in (\ref{new:5}). Making use of the estimate
\[
(\Delta_{\eta_t}^m{\cal R}_{\bar\alpha})(y)\leq
c\rho^{\alpha^*-n-m\frac{\alpha^*}{\max\limits_j\alpha_j}}(y),\
\rho(y)>A,\
\eta_t=t/\rho^\lambda(t)
\]
(see \cite{b4}), we arrive at the relation  $(\Delta_t^m{\cal R}_{\bar\alpha})(y)\in
L_1$. Therefore the left-hand side of (\ref{new:5}) converges in
$L_{\bar p}$, as $\varepsilon \to 0$, to $T^{\overline{\alpha}}f$, while
the right-hand side
%of (\ref{new:5})
converges in $L_{\bar r}$ to
$\Delta_t^m f$, since $f(x)\in L_{\bar r}$ and ${\cal
K}(y)\in L_1$. This implies
(\ref{5.21}).\quad $\hfill\blacksquare$

\subsection{Proof of inequality (\ref{5.22}).}
% Corollary~\ref{l5.2}}

With the aid of (\ref{5.21}) we have
\[
\dst \|\Delta_t^m f\|_{\bar p}\leq
\rho^{\alpha^*}(t) \|T^{\bar\alpha}f\|_{\bar p}\int\limits_{{\Bbb
R}^n} |(\Delta_{\eta_t}^m{\cal R}_{\bar\alpha})(y)|dy,
\]
where
$\eta_t=t/\rho^\lambda(t)$.
%by the representation (\ref{5.21}) and the Minkowsky inequality.
Since the integral on the right-hand side
is finite, we obtain (\ref{5.22}).\quad $\hfill\blacksquare$

\subsection{Proof of (\ref{4.4})}

The application of (\ref{1.5}) yields
\[
 \left|\left(B_N f\right)(x) \right| \leq \frac{C}{N^{\theta}}
 \int\limits_{{\Bbb R}^n} \frac{\rho^{\theta}(t) |f(x-t)|\, dt}{\left( 1+
 \frac{\rho(x)}{N}\right)^{\theta} \left( 1+
 \frac{\rho(x-t)}{N}\right)^{\theta} \rho^{n+\alpha^{\ast}}(t)},
\]
where $\dst\theta=\frac{\alpha^*}{\max\limits_j \alpha_j}$. Applying
the Minkowsky and H\"older inequalities, we obtain
\begin{eqnarray*}
 &&  \left\|B_N f\right\|_{\overline{p}} \leq C
 \|f\|_{\overline{q}} \int\limits_{{\Bbb R}^n}
 \frac{dt}{\rho^{n+\alpha^{\ast}-\theta}(t)} \times\\
 &&\quad \times \left( \int\limits_{{\Bbb R}^1}\!\dots\! \left( \int\limits_{{\Bbb
 R}^1}\!
 \left(\,\int\limits_{{\Bbb R}^1}\! \frac{ dx_1}{(1+\rho(x))^{\theta
 \beta_1} (1+\rho(x-t))^{\theta \beta_1} }
 \right)^{\frac{\beta_2}{\beta_1}}\!\! dx_2
 \right)^{\frac{\beta_3}{\beta_2}}\!\!\dots dx_n
 \right)^{\frac{1}{\beta_n}},
\end{eqnarray*}
where $\dst \beta_i=\frac{p_iq_i}{q_i-p_i}$, $i=1,\dots,n$.
Since $\displaystyle \sum\limits_{i=1}^n \frac{1}{\alpha_i
\beta_i} < \frac{\theta}{\alpha^{\ast}}$ ,   we arrive at
(\ref{4.4}) by virtue of (\ref{1.2}). \quad $\hfill\blacksquare$

{\bf Acknowledgments.} We are thankful to Prof. S.\,G.~Samko for
the useful discussions of the results and to Prof. E. Liflyand for
many helpful comments to the manuscript of the paper.

\vspace{\baselineskip}
%This work has been supported by Russian Fund Fundamental
%Investigations under grant no.~00--01--00046a.

\noindent
e-mail: vnogin@math.rsu.ru \\
e-mail: ournyce@macs.biu.ac.il

\end{document}